
%
%
\input mssymb
\def\CC {{\Bbb C}}
\def\FF {{\Bbb F}}
\def\HH {{\Bbb H}}
\def\NN {{\Bbb N}}

\def\QQ {{\Bbb Q}}
\def\RR {{\Bbb R}}
\def\ZZ {{\Bbb Z}}
\def\contin {\subseteq}

\def\part#1#2 {{\partial {#1}/\partial {#2}}}

\def\Sing {\mathop{\rm Sing}\nolimits}

\def\NS {\mathop{\rm NS}\nolimits}

\def\Hom {\mathop{\rm Hom}\nolimits}

\def\Sp {\mathop{\rm Sp}\nolimits}

\def\im {\mathop{\rm Im}\nolimits}
\def\orb {\mathop{\rm orb}\nolimits}

\def\bdn {{\bf n}}

\def\bdw {{\bf w}}
\def\bdx {{\bf x}}

\def\bdz {{\bf z}}

\def\To{\longrightarrow}

\def\Sum{\sum\limits}

\def\tens{\otimes}

\def\pf{\noindent{\sl Proof: }}

\font\gothic=eufm10
\font\head=cmr12

\def \goths {{\hbox{\gothic S}}}
\def \hf {{{1}\over{2}}}
\def\rationalmap{\mathrel{{\hbox{\kern2pt\vrule height2.45pt depth-2.15pt
 width2pt}\kern1pt {\vrule height2.45pt depth-2.15pt width2pt}
  \kern1pt{\vrule height2.45pt depth-2.15pt width1.7pt\kern-1.7pt}
  {\raise1.4pt\hbox{$\scriptscriptstyle\succ$}}\kern1pt}}}
\def\qed{\vrule width5pt height5pt depth0pt\par\smallskip}
\def\surj{\to\kern-8pt\to}
\outer\def\startsection#1\par{\vskip0pt
 plus.3\vsize\penalty-100\vskip0pt
  plus-.3\vsize\bigskip\vskip\parskip\message{#1}
   \leftline{\bf#1}\nobreak\smallskip\noindent}
\outer\def\camstartsection#1\par{\vskip0pt
 plus.3\vsize\penalty-150\vskip0pt
  plus-.3\vsize\bigskip\vskip\parskip
   \centerline{\it#1}\nobreak\smallskip\noindent}
\def\chain{\dot{\hbox{\kern0.3em}}}
\def\cochain{\d{\hbox{\kern0.3em}}}

\def\imic{\cong}
\outer\def\thm #1 #2\par{\medbreak
  \noindent{\bf Theorem~#1.\enspace}{\sl#2}\par
   \ifdim\lastskip<\medskipamount \removelastskip\penalty55\medskip\fi}
\outer\def\prop #1 #2\par{\medbreak
  \noindent{\bf Proposition~#1.\enspace}{\sl#2}\par
   \ifdim\lastskip<\medskipamount \removelastskip\penalty55\medskip\fi}
\outer\def\lemma #1 #2\par{\medbreak
  \noindent{\bf Lemma~#1.\enspace}{\sl#2}\par
   \ifdim\lastskip<\medskipamount \removelastskip\penalty55\medskip\fi}
\outer\def\corollary #1 #2\par{\medbreak
  \noindent{\bf Corollary~#1.\enspace}{\sl#2}\par
   \ifdim\lastskip<\medskipamount \removelastskip\penalty55\medskip\fi}
\def\deep #1 {_{\lower5pt\hbox{$\scriptstyle{#1}$}}}
%
\hyphenation{math-e-ma-ti-schen Wis-sen-scha-ft-en para-met-rised}
\def\restr #1 {{\big\vert_{#1}}}
\def\Ap {{{\cal A}_p}}
\def\Aplev {{{\cal A}_p^{\rm (lev)}}}
\def\PSL {\mathop{\rm PSL}\nolimits}
\def\At {{{\cal A}_p^*}}
\def\Atlev {{{\cal A}_p^{{\rm (lev)}*}}}
\def\cu#1#2 {{{\cal M}^*_#1(#2)}}
\def\Gp {{\Gamma_p}}
\def\Gplev{{\Gamma_p^{\rm (lev)}}}
\def\Trace {\mathop{\rm Trace}\nolimits}
\def\Fix {\mathop{\rm Fix}\nolimits}
\def\SL {\mathop{\rm SL}\nolimits}
\def\GL {\mathop{\rm GL}\nolimits}
\def\diag {\mathop{\rm diag}\nolimits}
\def\Stab {\mathop{\rm Stab}\nolimits}
\def\Vol {\mathop{\rm Vol}\nolimits}
\def\one {{\bf 1}}
\def\cH {{\cal H}}

\def\cC {{\cal C}}
\def\cF {{\cal F}}
\def\cV {{\cal V}}
\def\Alev {{\bf A}^{\rm (lev)}}
\def\ch {\mathop{\rm ch}\nolimits}
\def\Ahat {{\hat{\cal A}_p}}
\def\ratmap{\relbar\mkern3mu\to}
\def\split {{\mathrel>\joinrel\mathrel\triangleleft}}
\raggedbottom
\baselineskip=24pt

\raggedbottom
\magnification=1200

\centerline{\head MODULI OF POLARISED ABELIAN SURFACES}
\medskip
\centerline{\it G.K. Sankaran, Cambridge}
\medskip
\noindent Abelian surfaces over $\CC$ with a polarisation of type
$(1,t)$, $t$~a positive integer, are parametrised by a coarse moduli
space ${\cal A}_t$ which is a quasiprojective variety. In this paper
we shall concentrate on the case where $t$ is a prime $p\ge 5$, and
show that for $p$ sufficiently large (in fact $p\ge 173$) any
algebraic compactification of $\Ap$ is of general type.

A few results similar to this are already known. O'Grady, in [O'G],
considers the case $t=p^2$ and shows that a compactification of
${\cal A}_{p^2}$ is of general type for $p\ge 17$ (improved to $p\ge 11$
in~[GS]). The special feature here is the existence of a finite
morphism from ${\cal A}_{p^2}$ to the moduli space of principally polarised
abelian surfaces (the case $t=1$). This is also the case in [Bor],
where it is shown that a compactification of a Siegel modular
threefold coming from a subgroup~$\Gamma<\Sp(4,\ZZ)$ of finite index is of%

general type except for finitely many~$\Gamma$.
Another moduli space, referred to in this paper as $\Aplev$ and
parametrising abelian surfaces with a polarisation of type $(1,p)$ and
a level structure, is studied in depth by Hulek, Kahn and Weintraub in
the book [HKW2]. Its singularities are described in [HKW1] and it has
been shown, by Hulek, Gritsenko and me, that it is of general type if
the prime~$p$ is at least~$37$: see [HS] and [GH].

There is a finite morphism $\Aplev\to\Ap$. This morphism, the
singularities of $\Ap$ and its toroidal compactifications have been
studied in [Br] by Brasch, who gives an analysis in
the spirit of [HKW2]. Our main tools are Brasch's results, the
calculations relating to $\Aplev$ found in [HKW1], [HKW2] and [HS], and
some special cusp forms constructed by Gritsenko (see [G]). In
principle we do not need to know about $\Aplev$ but, like Brasch, we
are able to save a considerable amount of effort by making use of the
existing knowledge of it.

The paper is organised as follows. In Section~1 we give a more
detailed outline of the proof, collecting necessary facts from
elsewhere and establishing notation. Section~2 is devoted to
estimating the dimension of the space of cusp forms of high weight for
the relevant subgroup of $\Sp(4,\QQ)$. In Section~3 we describe a
method, due to Gritsenko, that enables us to handle the smooth points
at infinity in a suitable compactification of $\Ap$. Section~4
deals with the obstructions arising inside $\Ap$. Section~5 deals
with the singularities at infinity, and in Section~6 we assemble all
the parts of the proof and add some other remarks and corollaries.

I should like to thank Klaus Hulek, Valeri Gritsenko and J\"org Zintl
for helpful conversations, and the Deutsche Forschungsgemeinschaft and
the Max-Planck-Gesellschaft for financial support during my visits to
Hannover and Bonn in 1993 and 1994.

\startsection 1 Preliminaries.

Suppose $\Lambda\contin\CC^2$ is a lattice of rank~$4$. The complex
surface $S=\CC^2/\Lambda$ is called an abelian surface if it admits an
ample line bundle. A polarisation on an abelian surface $S$ is a class
$\lambda\in \NS(X)$ such that $\lambda=c_1({\cal L})$ for some ample
line bundle $\cal L$ on~$S$: see [LB] for details. As is well known
(see [LB] or [Mum]) a polarisation corresponds to an alternating
nondegenerate integral bilinear form on $\Lambda$, which can be
expressed, by choosing a suitable $\ZZ$-basis of $\Lambda$, by the
matrix
$$
\pmatrix{\phantom{-}0&T\cr -T&0\cr}
$$
where $T={\rm diag\ }(t_1,t_2)$, for some positive integers $t_1$
and $t_2$ with $t_1|t_2$. The integers $t_1$ and $t_2$ are uniquely
determined by $\lambda$, which is said to be a polarisation of type
$(t_1,t_2)$, or a $(t_1,t_2)$-polarisation. We may as well suppose
that $\lambda$ is not divisible in $\NS(S)$; then $t_1=1$ and $\lambda$
is a polarisation of type $(1,t)$ for some positive integer~$t$.

Denote by $\HH_2$ the Siegel upper half-plane of degree~$2$,
$$
\HH_2=\left\{ Z=\pmatrix{\tau_1&\tau_2\cr \tau_2&\tau_3\cr}\bigm|
Z\in M_{2\times 2}(\CC),\> Z=\!{}^\top\! Z,\> \im Z>0\right\},
$$
and define the paramodular group $\Gamma_t$ to be the subgroup of the
symplectic group $\Sp (4,\QQ)$ given by
$$
\Gamma_t=\left\{\gamma\in\Sp(4,\QQ)\Bigm|\gamma\in\pmatrix{
\phantom{t}\ZZ&\phantom{{1}\over{t}}\ZZ&\phantom{t}\ZZ&t\ZZ\cr
t\ZZ&\phantom{{1}\over{t}}\ZZ&t\ZZ&t\ZZ\cr
\phantom{t}\ZZ&\phantom{{1}\over{t}}\ZZ&\phantom{t}\ZZ&t\ZZ\cr
\phantom{t}\ZZ&{{1}\over{t}}\ZZ&\phantom{t}\ZZ&\phantom{t}\ZZ\cr}
\right\}.
$$
$\Gamma_t$ acts on $\HH_2$ by fractional linear transformations:
$$
\gamma=\pmatrix{A&B\cr C&D\cr}:Z\longmapsto (AZ+B)(CZ+D)^{-1}.
$$
The action is properly discontinuous and we denote the
quotient $\Gamma_t\backslash\HH_2$ by ${\cal A}_t$.

\prop 1.1 ${\cal A}_t$ is a coarse moduli space for pairs $(S,\lambda)$,
where $S$ is an abelian surface over~$\CC$ and $\lambda$~is a
polarisation of type~$(1,t)$.

\pf See [HKW2], Part~I, Chapter~1, where the prime~$p$ may be replaced
by an arbitrary positive integer with no other consequential changes.~\qed

Similarly we introduce the group
$$
\Gamma_t^{\rm(lev)}=\left\{\gamma\in\Sp(4,\QQ)\Bigm|
\gamma-\one_4\in\pmatrix{
\phantom{t}\ZZ&\phantom{t}\ZZ&\phantom{t}\ZZ&t\ZZ\cr
t\ZZ&t\ZZ&t\ZZ&t^2\ZZ\cr
\phantom{t}\ZZ&\phantom{t}\ZZ&\phantom{t}\ZZ&t\ZZ\cr
\phantom{t}\ZZ&\phantom{t}\ZZ&\phantom{t}\ZZ&t\ZZ\cr}
\right\}
$$
(in general we use $\one_n$, or just $\one$, to denote the $n\times n$
unit matrix) and the moduli
space ${\cal A}_t^{\rm (lev)}=\Gamma_t^{\rm (lev)}\backslash\HH_2$ of
$(1,t)$-polarised abelian surfaces with a level structure.

\prop 1.2 $\Gamma_t$ and $\Gamma_t^{\rm (lev)}$ are related as
follows:
\hfil\break
\qquad {\rm a)\ } $\Gamma_t^{\rm (lev)}$ is a normal subgroup of
$\Gamma_t$;
\hfil\break
\qquad {\rm b)\ } $\Gamma_t/\Gamma_t^{\rm (lev)}\imic \SL(2,\ZZ/t)$;
\hfil\break
\qquad {\rm c)\ } $\Gamma_t$ is conjugate in $\Sp(4,\QQ)$ to a
subgroup of $\Sp(4,\ZZ)=\Gamma_1$ if and only if $t$~is a square;
\hfil\break
\qquad {\rm d)\ } for an odd prime $p$ there is a morphism
$\Aplev\to\Ap$ of degree $p(p^2-1)/2$ exhibiting $\Ap$ as a
quotient of $\Aplev$ by an action of $\PSL(2,\ZZ/p)$.

\pf (a), (b) and (d) are in [HKW2] and (c) is easily checked.~\qed

{}From now on we shall concentrate on the case where $t$ is an odd prime
$p\ge 5$. Denote by $\At$ and $\Atlev$ the toroidal compactifications
of $\Ap$ and $\Aplev$ constructed in [Br] and [HKW2] respectively
(also called Igusa compactifications). Note that these constructions
are compatible with each other in the sense that the group action and
morphism of Propsition~1.2(d) extend to $\Atlev$ and to
$\Atlev\to\At$.

For any arithmetic subgroup $\Gamma$ of $\Sp (4,\QQ)$ we write $\cu
{k}{\Gamma} $ for the space of cusp forms of weight~$k$
for~$\Gamma$. Following Tai~([T]), we consider modular forms of
weight~$3n$ and look at $F\omega^{\tens n}$, where $F\in \cu
{{3n}}{\Gp} $ and $\omega=d\tau_1\wedge d\tau_2\wedge d\tau_3$ is a
differential $3$-form on~$\HH_2$. It is $\Gp$-invariant so it gives an
$n$-canonical form where the map $\HH_2\to\Ap$ is unbranched. To
obtain pluricanonical forms on a resolution of singularities of $\At$
we have to be able to extend over the smooth points of the boundary,
over the branch locus (including the singularities of $\Ap$) and over
the singularities in the boundary. The first and third of these may be
dealt with by choosing special cusp forms, partly in the manner of
[GH] and [GS]: for the second we have to estimate the growth with $p$
of the space of obstructions. These three types of obstacle will be
tackled in sections 3--5, but first we need to know how many cusp
forms there are.

\startsection 2 Cusp forms for $\Gp$.

In order to proceed as we intend we need a nontrivial cusp form of
weight~$2$ and a plentiful supply of cusp forms of large weight for $\Gp$.

\prop 2.1 There exists a nontrivial cusp form of weight~$2$ for $\Gp$
if~$p>71$.

\pf By [G], Theorem 3, a Jacobi cusp form of weight~$2$ and index~$p$
can be lifted to a weight~$2$ cusp form for $\Gp$, i.e.
$$
\dim\cu {2}{\Gp} \ge\dim{\goths}^J_{2,p}
$$
where ${\goths}^J_{k,t}$ is the space of Jacobi cusp forms of
weight~$k$ and index~$t$. But there is a formula for the latter
dimension (see [EZ] and [SZ]):
$$
\dim{\goths}^J_{2,p}= \Sum_{j=1}^p\{1+j\}_6-
\biggl\lfloor{{j^2}\over{4p}}\biggr\rfloor
$$
where
$$
\{m\}_6=\cases{\bigl\lfloor{{m}\over{6}}\bigr\rfloor& if $m \not\equiv
1$ mod $6$\cr
            \bigl\lfloor{{m}\over{6}}\bigr\rfloor-1& if $m \equiv 1$
mod $6$.\cr}
$$
It is easy to check that this number is positive if $p> 71$.~\qed

\noindent{\it Remark:\/} In fact this shows that
$\dim{\goths}^J_{2,t}>0$ for all $t>180$ and many smaller~$t$.
\prop 2.2 The space of cusp forms of weight~$k$ for~$\Gp$ satisfies
$$
\dim\cu {k}{\Gp} ={{p^2+1}\over{8640}}k^3+O(k^2)
$$
for any odd prime~$p$.

\pf The corresponding result for $\Gplev$ is given in [HS], Proposition~2.1,
where it is shown that
$$
\dim\cu {k}{\Gplev} ={{p(p^4-1)}\over{17280}}k^3+O(k^2)
$$

We proceed in the same way: writing $\bar\Gamma(1)=\Sp(4,\ZZ)/\pm\one$
and $\Gamma(l)$ for the principal congruence subgroup of level~$l$, we
have
$$
\dim\cu {k}{\Gamma(l)} ={{k^3}\over{8640}}\left[\bar\Gamma(1):\Gamma(l)\right]
$$
and if $p^2|l$ then $\Gamma(l)\contin\Gplev\contin\Gp$. Furthermore,
$\Gamma(l)$ is a normal subgroup of $\Gp$: denote the quotient by
$\Gp(l)$. Then
$$
\cu {k}{\Gp} =\cu {k}{\Gamma(l)} ^{\Gp(l)}
$$
and we can estimate the dimension, as in [T], by using the method of
Hirzebruch in [Hir]. We have (cf. [T], [GS])
$$
\eqalign{
\dim\cu {k}{\Gp} &=\dim\cu {k}{\Gamma(l)} ^{\Gp(l)}\cr
               &={{1}\over{|\Gp(l)|}}\sum_{\gamma\in\Gp(l)}
\Trace\Big(\gamma^*\restr {\cu {k}{\Gamma(l)} } \Big)
\cr}
$$
and by the Atiyah-Bott fixed point theorem
$$
\Trace\Big(\gamma^*\restr {\cu {k}{\Gamma(l)} } \Big)=O(k^{\dim\Fix(\gamma)})
$$
so we need to consider only $\gamma=\pm\one$, as otherwise
$\dim\Fix(\gamma)<3$.

But $-\one$ acts trivially, so we get
$$
\eqalign{
\dim\cu {k}{\Gp} &={{2}\over{|\Gp(l)|}}\dim\cu {k}{\Gamma(l)} \cr
               &={{2}\over{\left[\Gp:\Gamma(l)\right]}}{{k^3}\over{8640}}
 \bigl[\bar\Gamma(1):\Gamma(l)\bigr] + O(k^2)\cr
               &={{\bigl[\bar\Gamma(1):\Gplev\bigr]}\over
{\bigl[\Gp:\Gplev\bigr]}}{{k^3}\over{4320}} + O(k^2)\cr
               &={{p(p^4-1)/2}\over{|\SL(2,\ZZ/p)|}}{{k^3}\over{4320}}
+O(k^2)\cr
               &={{p^2+1}\over{8640}}k^3 + O(k^2)\cr}
$$
since $[\bar\Gamma(1):\Gplev]=p(p^4-1)/2$ by [HW], p.413,
and $\Gp/\Gplev\imic\SL(2,\ZZ/p)$ (Proposition~1.2), which has order
$p(p^2-1)$.~\qed

\noindent{\it Remarks.\/} a) The degree of the covering $\Aplev\to\Ap$ is
$p(p^2-1)/2$ because $-\one\in\Gp$ but $-\one\not\in\Gplev$. Thus the Galois
group is $\PSL(2,\ZZ/p)$, not $\SL(2,\ZZ/p)$.

\noindent\phantom{{\it Remarks.\/}} b) This is the first place where
we have assumed that $p$ is prime.Compare the corresponding
calculation for $t=p^2$ in~[GS].

\startsection 3 Extension of differential forms over the smooth part of
the boundary

Suppose $F_2$ is a nontrivial weight~$2$ cusp form for $\Gp$ and $F_n$
is a cusp form of weight~$n$. Then $F=F_2^nF_n$ is a cusp form of
weight~$3n$. The differential form $F\omega^{\tens n}$ on $\HH_2$
descends under the action of $\Gp$ to give an $n$-canonical form on
$\Ap$ away from the branch locus of $\HH_2\to\Ap$. Because $F$ has
been chosen to vanish to high order at infinity we get even more.

\prop 3.1 The differential form on $\Ap$ coming from $F\omega^{\tens n}$
extends holomorphically over the generic point of each codimension~$1$
boundary component of~$\At$.

\pf This is a straightforward application of [SC], Chapter~IV,
Theorem~1. The details are in [GS] (Proposition~3.2) for the case
$t=p^2$, and as the case $t=p$ is no different we omit them here.~\qed

Note that we do not need to know anything about the codimension~$1$
boundary components, not even how many there are (two, in fact). Of
course the forms also extend over smooth points in boundary components
of codimension greater than ~$1$, so the remaining obstructions to
extension come from the branch locus of $\HH_2\to\Ap$ and from the
singularities of $\At$.

\startsection 4 Branch locus and singularities of $\Ap$.

The singularities of $\At$ and the branch locus of $\HH_2\to\Ap$ are
fully described in [Br]. In this section we shall not deal with the
boundary $\At\setminus\Ap$, except when we deal with the points that
lie in the closure of the branch locus. Components of the singular
locus of $\At$ which are entirely contained in the boundary will be
dealt with in the next section.

We want to make use of the information about $\Atlev$ which is
already available in [HKW1], [HKW2] and [HS] (and also [Z], though we
shall not need that information directly). For this reason we need a
morphism $\Atlev\to\At$ which we can control, and which extends
the quotient morphism $\Aplev\to\Ap$. Perhaps there is an Igusa
compactification functor sending an arithmetic subgroup of
$\Sp(4,\QQ)$ to the Igusa compactification of the corresponding Siegel
modular $3$-fold and respecting inclusions among such groups. In
default of a result to that effect we prove the following, which is
sufficient for our immediate purposes.

\prop 4.1 The group $G=\Gp/\Gplev\imic\SL(2,\ZZ/p)$ acts on the Igusa
compactification $\Atlev$. The action agrees with the natural action
of $G$ on $\Aplev$ and the quotient is isomorphic to the Igusa
compactification~$\At$ of~$\Ap$.

\pf The Igusa compactifications $\Atlev$ and $\At$ are by
definition the compactifications constructed in [HKW2] (Part~I,
Chapter~3) and in [Br] respectively. The existence of a $G$-action on
$\Atlev$ extending the action on $\Aplev$ is an immediate
consequence of the construction of $\Atlev$: see Remark~3.95 on
page~89 of [HKW2]. So we must check that the quotient agrees with
Brasch's construction. By [HKW2], Part~I, Proposition~3.154, the
quotient $\Atlev/G$ is the compactification of $\Ap$ determined by
the Legendre decomposition. But this is precisely what is constructed
in Chapter~3,~{\S 2} of~[Br].~\qed

Both $\Gp$ and $\Gplev$ contain elements of order~$2$, among others
$$
I_1=\pmatrix{
          -1&0&\phantom{-}0&0\cr
\phantom{-}0&1&\phantom{-}0&0\cr
\phantom{-}0&0&          -1&0\cr
\phantom{-}0&0&\phantom{-}0&1\cr
},\qquad\qquad
I_2=\pmatrix{
          -1&          -1&\phantom{-}0&\phantom{-}0\cr
\phantom{-}0&\phantom{-}1&\phantom{-}0&\phantom{-}0\cr
\phantom{-}0&\phantom{-}0&          -1&          -1\cr
\phantom{-}0&\phantom{-}0&\phantom{-}0&\phantom{-}1\cr
}
$$
which are in $\Gplev$ and hence in $\Gp$, and $-\one_4$ which is in
$\Gp$~only. (Here we are using the notation and conventions of [HKW1],
which differ slightly from those of~[Br].)

We put
$$
\eqalign{
\cH_1&=\Fix(I_1)=\left\{\pmatrix{
\tau_1&0\cr
0&\tau_3\cr
}
\in\HH_2\mid\tau_1,\tau_3\in\HH_1\right\},\cr
\cH_2&=\Fix(I_2)=\left\{\pmatrix{
\phantom{-\hf}\tau_1&         -\hf\tau_3\cr
          -\hf\tau_3&\phantom{-\hf}\tau_3\cr
}
\in\HH_2\mid\tau_1,\tau_3\in\HH_1\right\}.\cr
}
$$

Denote by $H_i$ and $ H_i^{\rm (lev)}$ for $i=1,2$ the closures of the
images in $\cH_i$ in $\At$ and $\Atlev$ respectively. These are
the Humbert surfaces referred to in [HKW1], [HKW2] and [Br]. $\cH_1$
parametrises products of elliptic curves and $\cH_2$ parametrises
bielliptic abelian surfaces.

\prop 4.2 Every involution (element of order~$2$ different from
$-\one_4$) in $\Gp$ is conjugate in $\Gp$ to $\pm I_1$ or $\pm I_2$.
Every involution in $\Gplev$ is conjugate to $I_1$ or~$I_2$.

\pf See [Br] and [HKW1].~\qed

Recall that $-\one_4$ acts trivially on $\HH_2$, so that $\Gp$ acts on
$\HH_2$ through the effective action of the quotient $\bar\Gamma_p =
\Gp/\pm\one_4$. The elements $I_1$ and $I_2$ act by
reflection near $\cH_1$, $\cH_2$, so in view of Proposition~4.2 the
maps $\HH_2\to\Ap$ and $\HH_2\to\Aplev$ are branched over $H_i$
(respectively $ H_i^{\rm (lev)}$) and the singular locus, but nowhere else.

$\At$ and $\Atlev$ are normal projective varieties with only
finite quotient singularities. In particular they are smooth in
codimension~$1$, so $H_i$ and $ H_i^{\rm (lev)}$ are the only branch
divisors, and are $\QQ$-Cartier Weil divisors.

\corollary 4.3 If $n\in\NN$ is sufficiently divisible, $F_2\in \cu
{2}{\Gp} $ and $F_n\in \cu {n}{\Gp} $, then $n(K_{\At}+\hf H_1 +\hf
H_2)$ is Cartier and $F_2^nF_n\omega^{\tens n}$ determines an element
of $H^0\big(\At,n(K_{\At}+\hf{H_1}+\hf{H_2})\big)$.

\pf Clear: cf. [HS], Proposition 4.2.~\qed

\corollary 4.4 We may also also consider $F_2$ and $F_n$ as cusp forms for
$\Gplev$ and $F_2^nF_n\omega^{\tens n}$ as a $G$-invariant form with
poles on $\Aplev$, namely
$$
F_2^nF_n\omega^{\tens n} \in H^0\big(\At,
n(K_{\At}+\hf H^{\rm (lev)}_1+\hf H^{\rm (lev)}_2)\big)^G.
$$

\pf Obvious.~\qed

As well as the Humbert surfaces, we must consider the non-canonical
singularities of $\At$, which also provide an obstruction to
extending our pluricanonical forms to the whole of a smooth model.

\prop 4.5 The non-canonical singularities of $\At$ lie either in the
boundary of $\At\setminus\Ap$ or in~$H_1$.

\pf [Br], Hauptsatz.~\qed

For the remainder of this section we shall be concerned with what
happens in~$H_1$.

\prop 4.6 The non-canonical singularities of $\Ap$ lying in $H_1$ are
precisely the points of of the two curves $C_{3,1}$ and $C_{5,1}$,
which are the images in $\Ap$ of
$$
\eqalign{
\cC_{3,1}&=\left\{\pmatrix
{\tau_1&0   \cr
      0&\rho\cr}
\Bigm|\tau_1\in\HH\right\}\cr
\cC_{5,1}&=\left\{\pmatrix
{\rho&0     \cr
    0&\tau_3\cr}
\Bigm|\tau_3\in\HH\right\}\cr}
$$
where $\rho=e^{2\pi i/3}$. The transverse singularity at the generic
point of each of these curves is the cone on the twisted cubic curve.

\pf According to [Br], Hilfs\"atze 2.24 and 2.25, $H_1\cap\Sing\Ap$
consists of four curves, $C_{2,1}$, $C_{3,1}$, $C_{4,1}$ and
$C_{5,1}$, forming a square (see Figure~1).
\midinsert
$$\vbox{\offinterlineskip
\halign{
#\hskip 1.5cm&\vrule#&#\hskip 1cm&\hfil#\hfil&#\hskip 1
cm&\vrule#&#\hskip 1.5cm\cr
\strut&&&&&&\cr
\strut&&&&&&\cr
&&$\,P_{1,1}$&&&$\,P_{4,1}$\cr
\noalign{\hrule}
&height 2pt&&&&&\cr
&&&$C_{3,1}$&&&\cr
\strut&&&&&&\cr
\strut&&&&&&\cr
\phantom{$\,C_{5,1}$}&&$\,C_{5,1}$&\phantom{$\,C_{6,1}$}&&$\,C_{6,1}$\cr
\strut&&&&&&\cr
\strut&&&&&&\cr
&&&$C_{4,1}$&&&\cr
&height 2pt&&&&&\cr
\noalign{\hrule}
&height 2pt&&&&&\cr
&&$\,P_{3,1}$&&&$\,P_{2,1}$\cr
\strut&&&&&&\cr
\strut&&&&&&\cr
\multispan7\hfil\strut\cr
\multispan7 \hfil Figure~1\hfil\cr}}
$$
\endinsert

The transverse singularity at a general point of $C_{4,1}$ or
$C_{6,1}$ (specifically, at a point of either of these curves which is
not also in another $C_{i,1}$ nor in the boundary of $\At$) is an
ordinary double point. The isotropy group $Z(x)$ of such a point is of
order~$8$ but includes $-\one_4$ and $I_1$, which generate a normal
subgroup of index~$2$. At $x$, $-\one_4$ acts trivially on the tangent
space and $I_1$ acts by reflection, so the singularity is isomorphic
to a quotient singularity by an action of
$Z(x)/\langle-\one_4,I_1\rangle$, i.e., a threefold ordinary double
point. Such a singularity is, of course, canonical.

At a general point of $C_{3,1}$ or $C_{5,1}$ the isotropy group has
order~$12$ and again $-\one_4$ and $I_1$ generate a normal subgroup,
of index~$3$. We could calculate the transverse singularity directly,
but it is better to argue as follows. Suppose $x\in C_{5,1}$ but
$x\not=P_{3,1}$ or $P_{1,1}$ (i.e. $x$ is in no other component of
$\Sing\Ap$). Let $\tilde x \in \Aplev$ be such that $\tilde x\mapsto
x$ under $\Aplev \to \Ap$: then $\tilde x\in C_2$, in the notation of
[HKW2]. The isotropy group $Z(\tilde x)$ of $\tilde x$ in $\Gplev$ has
order~$6$ and does not contain $-\one_4$, so $Z(x)=\langle
Z(\tilde x),-\one_4\rangle$. But $-\one_4$ acts trivially so the singularity
at $x\in \Ap$ is the same as the singularity at $\tilde x \in\Aplev$,
and this is of transverse type the cone on the twisted cubic by [HKW2],
Theorem~1.8. If instead $x\in C_{3,1}$, we use the existence of an
extra automorphism $\Theta$ on $\Ap$ which sends an abelian variety to
its dual: it is represented by the element
$$
\Theta=\pmatrix{
0&{{1}\over{\surd p}}&0&0\cr
\surd p&0&0&0\cr
0&0&0&\surd p\cr
0&0&{{1}\over{\surd p}}&0\cr}
\in\Sp(4,\RR)
$$
acting on $\HH_2$ (one verifies that $\Theta$ normalises $\Gp$ in
$\Sp(4,\RR)$). This automorphism interchanges $C_{3,1}$ and $C_{5,1}$,
which therefore have the same singularities.

Finally, the singularity at the corner $P_{2,1}=C_{4,1}\cap C_{6,1}$
is canonical and the singularities at the other three corners are not,
by the criterion of Reid, Shepherd-Barron and Tai: see [YPG]~\S 4.11.~\qed

Let $\phi:\Atlev\to\At$ be the quotient morphism of
Proposition~4.1. Let $\beta_1^{\rm (lev)}:\Alev=
\tilde{\cal A}_p^{{\rm (lev)}*}\to\Atlev$ be the resolution of
singularities defined in [HS],~p.18: that is, blow up along $C_1$ and
$C_2$ in $\Atlev$ and take a $G$-invariant resolution of the other
singularities of $\Atlev$. Then there is a diagram
$$
\matrix{
\Alev=&\tilde{\cal A}_p^{{\rm (lev)}*}&
{\buildrel {\tilde\phi} \over {\relbar\joinrel\relbar\joinrel\To}}& \tilde{\cal
A}_p^{{\rm (lev)}}&\cr
&\scriptstyle{\beta_1^{\rm (lev)}}\Bigg\downarrow
\phantom{\scriptstyle{\beta_1^{\rm (lev)}}}&
&\phantom{{\scriptstyle{\beta_1}}}\Bigg\downarrow
\scriptstyle{\beta_1}&\cr
&\Atlev&{\buildrel \phi \over {\relbar\joinrel\relbar\joinrel\To}}&\At&\cr}
$$
where $\phi$ and $\tilde\phi$ are quotient maps under the action of
$G$, and $\beta_1$ factors through the blow-up of $\At$ along
$\phi(C_1)=C_{6,1}$ and $\phi(C_2)=C_{5,1}$. ($\beta_1$ also makes
some other modifications in the boundary, which do not concern us.)
$\Alev$ is the same as the smooth variety $\bf A$ studied in \S 4
of~[HS]. The singularities at the general points of $C_{5,1}$ and
$C_{6,1}$ are resolved by~$\beta_1$.

We have used $C_{5,1}$ (etc.) to denote both a curve in $\Ap$ and its
closure in $\At$. To this abuse of notation we add another: we will
denote the strict transforms of $H_1$ and $ H_i^{\rm (lev)}$ under $\beta_1$
and $\beta_1^{\rm (lev)}$ by $H_i$ and $ H_i^{\rm (lev)}$ again. We denote the
exceptional divisor of $\beta_1$ over $C_{5,1}$ by $E_{5,1}$, and so on, and
similarly the exceptional divisor of $\beta_1^{\rm (lev)}$ over $C_i$
will be $E_i$; thus $E_1$ and $E_2$ are what are called $E$ and $E'$
in [HS], and $\tilde\phi(E_1)=E_{6,1}$, $\tilde\phi(E_2)=E_{5,1}$.

\prop 4.7 The obstructions to extending pluricanonical forms over
the general points of $E_{6,1}$, $H_1$ and $H_2$ are zero.

\pf The singularities along $C_{6,1}$ are canonical at the general
point. The obstructions coming from $E_1$, $H^{\rm (lev)}_1$ and
$H^{\rm (lev)}_2$ are in any case shown to be zero in [HS],
Corollary~4.7, Theorem~4.19 and Theorem~4.25. The obstructions we are
interested in are just the $G$-invariant parts of those obstructions,
hence also zero.~\qed

If $F_2$ and $F_n$ are cusp forms for $\Gp$ of weights~$2$ and~$n$
($n$ sufficiently divisible) then $F_2^nF_n\omega^{\tens n}$ gives,
exactly as in [HS], a $G$-invariant element
$$
\eta\in H^0\big(\Alev;n(K_{\Alev}+\hf H^{\rm (lev)}_1+\hf H^{\rm (lev)}_2
+{{1}\over{4}}E_1+\hf E_2)\big)^G
$$
and by Proposition 4.7, above, and Proposition 4.8 of [HS], we have at
once the following.

\prop 4.8 The form $\eta$ is in fact a $G$-invariant element
$$
\eta\in H^0(\Alev;nK_\Alev)^G
$$
provided it lies outside a subspace of dimension at most
$$
\sum_{j=1}^{n/2}\dim H^0\big(\bigl[n(K_{\Alev}+\hf H^{\rm (lev)}_1+\hf
H^{\rm (lev)}_2
+{{1}\over{4}}E_1+\hf E_2)-({{n}\over{2}}-j)E_2\bigr]\restr {E_2} \bigr)^G.
$$

\pf As for [HS], Proposition 4.3, but taking $G$-invariant
sections.~\qed

Put $n=12n'$ and $L_j=\bigl[n(K_{\Alev}+\hf H^{\rm (lev)}_1+\hf H^{\rm
(lev)}_2 +{{1}\over{4}}E_1+\hf
E_2)-({{n}\over{2}}-j)E_2\bigr]\restr {E_2} $. Let $\Sigma'$ and $\Phi'$ be
a section and a fibre of the ruled surface $E_2=E'$, as in~[HS].

\thm 4.9 The obstruction coming from $E_2$ to extending
$F_2^nF_n\omega^{\tens n}$, where $n$ is sufficiently divisible, is
$$
\left({{7}\over{108}}-{{1}\over{6p}}\right)n^3+O(n^2).
$$

\pf We want to calculate $\sum_{j=1}^{6n'}\dim H^0(L_j)^G$. By [HS],
Proposition~4.11,
$$
L_j\equiv(12n'-3j)\Sigma'+\bigl[6n'(\mu-2\nu_\infty)-j\mu/2\bigr]\Phi',
$$
where $\nu_\infty=(p^2-1)/12$ and $\mu=p\nu_\infty$: as in the proof of
[HS], Theorem 4.12, $H^0(L_j)=0$ for $j>4n'$ and $L_j-K_{E_2}$ is
ample for~$j\le 4n'$.

To estimate $\dim H^0(L_j)^G$ for $j\le 4n'$ we use the Atiyah-Bott
fixed point theorem, much as we did in \S 2 above. By an elementary
result about finite group representations (e.g. [Se],~I.2.3)
$$
\eqalign{
\dim H^0(L_j)^G
  &= {{1}\over{|G|}}\sum_{g\in G}\Trace\,(g^*\restr {{H^0(L_j)}} )\cr
  &= {{1}\over{|G|}}\sum_{g\in G}\sum_i(-1)^i\Trace\,(g^*\restr
{H^i(L_j)} )\cr
}
$$
since the higher cohomology vanishes. But by the fixed point theorem
([AS], [Hir]; cf.~[T], Appendix to \S 2)
$$
\sum_{g\in G}\sum_i(-1)^i\Trace\,(g^*\restr {H^i(L_j)} )
 =\Bigl\{\ch(L_j\restr {E_2^g} )(g)\cdot\mu_g\Bigr\}\bigl[E_2^g\bigr]
$$
where $\mu_g$ is a class depending only on $g$, not~$L_j$. This is a
polynomial whose total degree in $n'$ and~$j$ is the degree of
$\ch(L_j\restr {E_2^g} )(g)$ in $n'$ and~$j$, which is at most $\dim E_2^g$.
Therefore the only $g\in G$ which make a contribution involving
$n'^2$, $n'j$ or $j^2$ are~$\pm 1$. So
$$
\eqalign{
\sum_{j=1}^{6n'}\dim H^0(L_j)^G
  &= \sum_{j=1}^{4n'} {{1}\over{|G|}}\sum_{g=\pm 1}\sum_i(-1)^i
       \Trace\,(g^*\restr {H^i(L_j)} ) + O(n^2)\cr
  &= \sum_{j=1}^{4n'} {{1}\over{|G|}} \Big(\Trace\,(1^*\restr
      {H^0(L_j)} )+\Trace((-1)^*\restr {H^0(L_j)} )\Big) + O(n^2)\cr
  &= {{2}\over{|G|}}\sum_{j=1}^{4n'}\dim H^0(L_j) + O(n^2)\cr
  &= {{2}\over{|G|}}\Big({{7}\over{108}}\mu-{{1}\over{6}}\nu_\infty\Big)n^3
       + O(n^2)\cr
  &= \Big({{7}\over{108}}-{{1}\over{6p}}\Big)n^3 +O(n^2)\cr}
$$
using the calculation of $\dim H^0(L_j)$ in~[HS].~\qed

The obstruction coming from $C_{3,1}$ could be calculated directly in a
similar way but there is no need to do this. Instead we make use of
the extra symmetry~$\Theta$ that corresponds to interchanging an
abelian surface and its dual.

\thm 4.10 The obstruction to extending a form $\eta$ to a
pluricanonical form on a resolution of the singularities of $\At$
that lie in $\Ap$ is contained in a space of dimension at most
$$
\Big({{7}\over{54}}-{{1}\over{3p}}\Big)n^3 + O(n^2).
$$

\pf Consider the blow-up $\beta_2:\bar{\At}\to\tilde\At$ along
the curves $C_{3,1},C_{4,1}\contin\tilde\At$. Clearly this resolves
the singularities at the general points of $C_{3,1}$ and $C_{4,1}$. It
is easy to see that the singularity of $\bar{\At}$ above $P_{1,1}$ is
an ordinary double point and thus canonical: hence all the
singularities of $\bar{\At}$ away from the boundary are canonical.
But we could also blow up $C_{3,1}$ and $C_{4,1}$ first and then
$C_{5,1}$ and $C_{6,1}$. Let $\beta'_2:\tilde\At'\to\At$ and
$\beta'_1:\bar{\At}'\to\tilde\At'$ be these blow-ups. The left and
right halves of the diagram
$$
\matrix{
\bar{\At}&&&&\bar{\At}'\cr
&\searrow^{\beta_1\beta_2}&&{}^{\beta'_1\beta'_2}\swarrow&\cr
&&\At&&\cr}
$$
are interchanged by $\Theta$, so if $\cV$ is the space of forms~$\eta$
extending over a resolution of $C_{5,1}$ then $\Theta^*\cV$ is the
space of forms extending over a resolution of $C_{3,1}$. Furthermore,
if $\eta\in\cV\cap\Theta^*\cV$ then $\eta$ extends over the smooth
part of the open subsets of $\bar{\At}$ and $\bar{\At}'$ where the
birational map $\bar{\At}\ratmap\bar{\At}'$ is an isomorphism. But
this birational map is an isomorphism except over the corners
$P_{1,1}$, $P_{1,2}$, $P_{1,3}$ and $P_{1,4}$ in $\At$, and there
the singularities of $\bar{\At}$ are canonical, so such an $\eta$
extends as in the statement of the theorem. The codimension in
$\cu {n}{\Gp} $ of $\cV\cap\Theta^*\cV$ is at most twice the codimension
of $\cV$, which is what we need.~\qed

\noindent{\it Remarks.\/} a) Although we have not considered
singularities in the boundary, the blow-ups we have made do affect
them. We shall see later that they have changed them for the better,
from our point of view.

\noindent{\phantom{\it Remarks.\/}} b) In fact we do not expect the
conditions coming from $C_{3,1}$ and $C_{5,1}$ to be independent,
because the curves meet. But the space of common obstructions will be
small, that is to say $O(n^2)$, so we cannot gain anything by
calculating it.

\startsection 5 Singularities in the boundary.

Unlike $\Atlev$, the Igusa compactification $\At$ has
non-canonical singularities in the boundary, which present a further
obstacle to extending forms.

The boundary of $\Ap$ consists of two disjoint open surfaces,
$D^\circ(\ell_0)$ and $D^\circ(\ell_1)$ (equivalent under $\Theta$),
whose closures $D(\ell_0)$ and $D(\ell_1)$ meet in a corank~$2$
boundary component, a curve called~$E(h)$. This corresponds to the
structure of the Tits building of $\Gp$, shown in Figure~2. We shall
also use  $D(\ell_0)$, $D(\ell_1)$ and $E(h)$ to denote the strict
transforms of these subvarieties in $\bar{\At}$.
\midinsert
$$\displaylines{
\bullet\kern-0.1cm\hbox{\vrule width 2cm height 2.6pt depth -2pt}
\kern-0.1cm\bullet\kern-0.1cm\hbox{\vrule width 2cm height 2.6pt depth -2pt}
\kern-0.1cm\bullet\cr
\ell_0\hskip 1.8cm h\hskip 1.8cm \ell_1\cr
\hbox{Figure~2}\cr}
$$
\endinsert

\thm 5.1 All the non-canonical singularities of $\bar{\At}$ lie in
the corank~$2$ boundary component~$E(h)$.

\pf The singularities of $\At$ at the boundary are calculated in
[Br]. In $D^\circ(\ell_0)$ (respectively $D^\circ(\ell_1)$) there are
exactly four singular points, $Q_{1,0}$, $Q_{2,0}$, $Q_{3,0}$ and
$Q_{4,0}$ (respectively $Q_{1,1}$, $Q_{2,1}$, $Q_{3,1}$ and
$Q_{4,1}$). The singularities at $Q_{3,0}$, $Q_{4,0}$, $Q_{3,1}$ and
$Q_{4,1}$ are isolated, but the others are not: in fact
$Q_{1,0}=C_{4,1}\cap D^\circ(\ell_0)$, $Q_{2,0}=C_{3,1}\cap
D^\circ(\ell_0)$, $Q_{1,1}=C_{6,1}\cap D^\circ(\ell_1)$ and
$Q_{2,1}=C_{5,1}\cap D^\circ(\ell_1)$. All this is shown in [Br],
Kapitel~3, S\"atze~4.6,~4.7.

In the proofs of those theorems the isotropy groups are calculated. A
neighbourhood of $D^\circ(\ell_0)$ in $\At$ is isomorphic to a
neighbourhood of $t_1=0$ in $\CC\times\CC\times\HH/P''(\ell_0)$, where
$$
P''(\ell_0)=\left\{g''=\pmatrix{
\varepsilon&         \varepsilon m&         \varepsilon n\cr
          0&\phantom{\varepsilon}a&\phantom{\varepsilon}b\cr
          0&\phantom{\varepsilon}c&\phantom{\varepsilon}d\cr
}\Bigm|\varepsilon=\pm 1,\;(m,n)\in \ZZ^2E,\;
\pmatrix{a&b\cr c&d\cr}\in E^{-1}\SL(2,\ZZ)E\right\}
$$
for $E=\pmatrix{1&0\cr 0&p}\in\GL(2,\QQ)$. Here the action of $P''(h)$
on $(t_1,\tau_2,\tau_3)\in\CC\times\CC\times\HH$ is given by
$$
g'':\pmatrix{t_1\cr \tau_2\cr \tau_3\cr} \mapsto \pmatrix{
t_1\exp\bigl\{2\pi i[m\tau_2+(dm-c\tau_2-cn)\varepsilon\tau_3]\bigr\}
\cr
\varepsilon(\tau_2+m\tau_3+n)(c\tau_3+d)^{-1}
\cr
(a\tau_3+b)(c\tau_3+d)^{-1}
\cr
}.
$$
In particular $-\one_3$ acts trivially and we are really concerned
with the action of $P''(\ell_0)/\pm\one_3$.

$Q_{1,0}$ is represented by $(t_1,\tau_2,\tau_3)=(0,0,pi)$ and the
effective isotropy group is generated by
$$
\pmatrix{
1&0&0\cr 0&0&p\cr 0&-{{1}\over{p}}&0\cr}.
$$
This element has order~$2$ modulo $-\one_3$ and $I_1=\diag(-1,1,1)$
(which acts as a reflection and therefore introduces no singularity):
thus $Q_{1,0}$ is a quotient singularity of order~$2$, actually of
type $\hf(0,1,1)$ in the usual notation for cyclic quotient
singularities (see [YPG]). In any case it is resolved by~$\beta_2$.

$Q_{2,0}$ is represented by $(t_1,\tau_2,\tau_3)=(0,0,\rho)$, where
$\rho=e^{2\pi i/3}$ and the effective isotropy group is generated by
$$
\pmatrix{
1&0&0\cr 0&0&-p\cr 0&{{1}\over{p}}&1\cr}.
$$
This element has order~$3$ modulo $\left\langle-\one_3,I_1\right\rangle$
and the singularity is of type ${{1}\over{3}}(0,1,1)$, neither
isolated nor canonical. However, $Q_{2,0}$ is thus a point of
$C_{3,1}$ with the same singularity as the general point and is
therefore resolved by~$\beta_2$.

$Q_{3,0}$ is represented by
$(t_1,\tau_2,\tau_3)=(0,{{p}\over{2}}(i-1),pi)$,
and the isotropy comes from
$$
\pmatrix{
1&p&0\cr 0&0&p\cr 0&-{{1}\over{p}}&1\cr}.
$$
The singularity here is of type $\hf(1,1,1)$ (the cone on the
Veronese). One can see this by carrying out the same argument as in
[HKW1], Proposition~2.8 (for the point $Q'_1$), or by the following
direct algebraic argument. The point $Q_{3,1}$, which is of the same
type as $Q_{3,0}$, is the image of $Q'_1$ under
$\phi:\Atlev\to\At$. The orbit of $Q'_1$ under~$G$ consists of
$p^2-1$ points, one on each $D^\circ_{\ell(a,b)}$ (each peripheral
boundary component: see [HKW1]), so $|\Stab_G Q'_1|=2p$. Since
$-1\in\Stab_G Q'_1$ and acts trivially, the singularity at
$Q_{3,1}$ must either be the same as at $Q'_1$ or a quotient of it by
a cyclic group of order~$p$, depending on whether an element of
order~$p$ in $\Stab_G Q'_1$ has trivial image in $P''(\ell_0)$ or not
(i.e. whether it acts trivially on the Zariski tangent space or not).
But $P''(\ell_0)$ has no $p$-torsion.

Exactly the same argument shows that the singularity at $Q_{4,1}$, and
hence the one at $Q_{4,0}$, is the same as the singularity of
$\Atlev$ at $Q'_2$, which by [HKW1] is a cyclic quotient singularity of type
${{1}\over{3}}(1,2,1)$. Both $\hf(1,1,1)$ and ${{1}\over{3}}(1,2,1)$
are canonical singularities, so we have finished.~\qed

The singularities of $\At$ lying on $E(h)$ are described in [Br].
The picture is a little complicated, depending among other things on
the residue class of $p$ modulo~$12$. There are many (about $p/6$)
isolated cyclic quotient singularities present, including all the ones
of type ${{1}\over{p}}\big(r+1,-r,r(r+1)\big)$ where the residue class
of $r$ mod~$p$ is not $0$, $1$ or a primitive cube or fourth root of
unity in $\FF_p$. Such singularities are in general not canonical, as
follows from [Mor] and [MorS]. Moreover, the plurigenera of a
non-canonical 3-fold cyclic quotient singularity $P\in X$ of index~$p$
(that is, the dimension of the obstruction to extending sections of
$nK_X$ to a resolution) can be expected to be close to $p^2n^3$, and
thus almost as big as $\dim\cu {n}{\Gp} $. So a straightforward dimension
count is unlikely to succeed and we shall have to find a special
property of the pluricanonical forms we have chosen that allows them
to extend.

\noindent{\it Remark.}\/ To see why one expects the plurigenera to grow in
this way it is easiest to use toric methods, as described in [YPG].
Put $N=\ZZ^3$, $M=\Hom (N,\ZZ)$,
$N'=\ZZ^3+\ZZ\cdot{{1}\over{p}}(\nu_1,\nu_2,\nu_3)=\ZZ^3+\ZZ\cdot\bdn$ and
$M'=\Hom (N',\ZZ)$, with $0<\nu_i<p$. Then the plurigenus $P_n(\bdn)$
associated with the toric morphism corresponding to the ray $\RR\cdot\bdn$
is given by the number of points of $M'\cap n\Delta_{\bdn}^\circ$,
where
$$
\Delta_{\bdn}=\left\{\bdx=(\xi_1,\xi_2,\xi_3)\in\RR^3=M\tens\RR\mid
x_i>1, \langle\bdn,\bdx\rangle <1\right\}
$$
and $\Delta_{\bdn}^\circ$ denotes the Fine interior (see [YPG],
appendix to~{\S 4}). If we are only interested in the growth of
$P_n(\bdn)$ with~$n$ we may as well take the topological interior
instead. We have $|M:M'|=p$ and
$$
\Vol_M(\Delta_{\bdn})={{1}\over{6}}
\Bigl({{p}\over{\nu_1}}-1\Bigr)\Bigl({{p}\over{\nu_2}}-1\Bigr)
\Bigl({{p}\over{\nu_3}}-1\Bigr)
$$
so $P_n(\bdn)\sim {{p^2n^3}\over{6\nu_1\nu_2\nu_3}}$. If the
obstruction for arbitrary forms is not to dominate $\dim\cu {n}{\Gp} $,
therefore, we should need $\nu_1\nu_2\nu_3>1440$; and this is the
obstruction arising from just one blow-up on the way to resolving just
one of the singularities. In general it is clear that the obstruction
will be far too big.

We need a little more information about the singularities of $\At$ first.

\lemma 5.2 The non-canonical singularities of $\At$ lying in the
corank~$2$ boundary componenet $E(h)$ are all isolated and not in the
closure of the branch locus of the quotient map $\HH_2\to\Ap$.

\pf According to [Br], the singularities in $E(h)$ that lie in the
closure of the branch locus are either cyclic quotient singularities
of order~$2$, of type $\hf(1,2,1)$, or quotients by $\ZZ_p\split\ZZ_3$
in which generators $\zeta_3$ of $\ZZ_3$ and $\zeta_p$ of $\ZZ_p$ act by
$$
\eqalign{
\zeta_3:(z_1,z_2,z_3)&\mapsto(z_3,z_1,z_2)\cr
\zeta_p:(z_1,z_2,z_3)&\mapsto(e^{2\pi i(1+r)/p}z_1,e^{-2\pi
ir/p}z_2,e^{-2\pi i/p}z_3)\cr}
$$
for some $r\in\ZZ$. The first of these types occurs if $p\equiv 1$ or
$p\equiv 5$ mod~$12$, the second if $p\equiv 1$ or $p\equiv 7$ mod~$12$.
By the criterion of Reid, Shepherd-Barron and Tai these singularities
are canonical.~\qed

We are going to use the weight~$2$ form a second time. We assume, for
the moment, that we have a form that extends over the part we have
already covered and produce from it a form (of higher weight) that
extends everywhere.

\thm 5.3 Suppose that $n=3n'$ and that the form $F_n\omega^{\tens n'}$
extends to an $n'$-canonical form on a resolution of singularities of
$\At$ except perhaps over the exceptional set coming from the
singularities in the corank~$2$ boundary component $E(h)$ of $\At$.
Then the form $F_2^nF_n\omega^{\tens n}$ extends to an $n$-canonical
form on the resolution.

\pf Let $\psi:\Ahat^*\to\At$ be a resolution of the non-canonical
singularities in~$E(h)$. We may assume that the exceptional locus of $\psi$
is a normal crossings divisor. It is enough to show that
$F_2^nF_n\omega^{\tens n}$ extends to a pluricanonical section near a
general point of any exceptional component~$E$. Suppose, then, that
$P\in E$ and that $P$ is not in any other exceptional component. Let
$\alpha$ be the discrepancy at~$E$: that is,
$K_{\Ahat^*}-\psi^*K_{\At}=\alpha E+E'$, where $E'$ is supported on
the other exceptional components. If $\alpha\ge 0$ there is nothing to
prove, so assume $\alpha<0$. Then we can choose an analytic open
neighbourhood $U$ of $P$ in $\Ahat^*$ such that $U$ does not meet the
branch locus of $\pi:\HH_2\to \Ap$, by Lemma~5.2. Put $\tilde
U=\pi^{-1}(U)\cap\cF_p$, where $\cF_p$ is a fundamental domain for the
action of $\Gp$ on $\HH_2$ (there is an explicit description of
$\cF_p$ in [Br], following [Fr]). Then $\pi\restr {\tilde U} :\tilde U\to
U\setminus E$ is an isomorphism and we can identify $\tilde U$ with
$U\setminus E$.

We have $F_2^nF_n\omega^{\tens n}=(F_2^3\omega^{\tens 2})^{\tens
n'}\tens(F_n\omega^{\tens n'})$. Choose local
coordinates $(z_1,z_2,z_3)$
near $P$ so that $E$ is given by $z_1=0$ and take $U$ to be the
polycylinder $\{|z_1|<\epsilon,|z_2|<\delta,|z_3|<\delta\}$, so
${\tilde U}\imic\Delta_\epsilon^*\times\Delta_\delta^2$. If we take a
triple cover $\hat U\to U$ branched along $E$, by setting
$(w_1^3,w_2,w_3)=(z_1,z_2,z_3)$, we can write, formally
$$
F_2\omega^{2/3}=f(\bdw)(dz_1\wedge dz_2\wedge dz_3)
$$
(for $\hat U$ sufficiently small), and there is a Laurent series
$$
f(\bdw)=\sum_{r=-\infty}^{\infty} a_r w_1^r
$$
where the $a_r$ are analytic functions of $w_2$ and $w_3$. This is
because on $U$ we have $f=F_2J^{2/3}$, where $J=\det\left(\partial
\tau_i/ \partial z_j\right)$ is the Jacobian determinant. $J(\bdz)$ is
finite and nonzero on $\tilde U$ (since $\tau_i$ and $z_j$ are local
coordinates at each point) so we can take a cube root of $J$ on the
triple cover~$\hat U$.

I claim that $a_r=0$ for $r\le 0$, that is, that the order of
vanishing $v_E(F_2^3J^2)$ of $F_2^3J^2$ along $E$ is
positive. Certainly $v_E(F_2)>0$ as $F_2$ is a cusp form, so
$v_E(F_2^3)\ge 3$. We need to understand~$J$, which means having a
good description of the geometry near~$E$. So we need to look a little
more closely at the construction of the boundary.

Near $E(h)$ the structure of $\At$ is toroidal: it is the quotient
by a certain discrete group $P''(h)$ of an open subset of a torus
embedding $T\deep {N} {\rm emb}(\Sigma_h)$, where $\Sigma_h$ is a suitable
fan. The singularities that we are concerned with arise because the
fan $\Sigma_h$ that one naturally chooses is not basic (see [D] or
[Od] for the terminology of torus embeddings). There are also some
singularities arising from the fixed points of torsion elements in
$P''(h)$ but they are canonical. So by taking an equivariant
subdivision $\Sigma^+_h$ of $\Sigma_h$ we can obtain a resolution of
singularities of the type we want. Our exceptional component~$E$ then
comes from $\orb(\sigma)$ for some $1$-dimensional cone
$\sigma\in\Sigma^+_h\setminus\Sigma_h$. A local equation for $E$ at
$P$ is thus $t_1^{b_1}t_2^{b_2}t_3^{b_3}=0$, for suitable $b_i\in\ZZ$,
where the $t_i$ are coordinate functions on the torus
$T_N\imic(\CC^*)^3$. This will be valid on some affine open toric
variety, one of the pieces of $T_N{\rm emb}(\Sigma_h)$.

So we may take $z_1=t_1^{b_1}t_2^{b_2}t_3^{b_3}$. However, the torus
embedding is defined by setting $t_1=e^{2\pi i \tau_1}$, $t_2=e^{2\pi
i \tau_2/p}$ and $t_3=e^{2\pi i \tau_3/p}$, so we have
$$
z_1=e^{2\pi i\{b_1\tau_1+b_2\tau_2/p+b_3\tau_3/p\}}
$$
and $\partial z_1/\partial\tau_j=2\pi ib_jz_i$. Hence $J^{-2}$ is of
degree at most~$2$ in $z_1$, and thus $v_E(J^2)\ge -2$. Therefore
$v_E(F_2^3J^2)\ge 1$.

What this means is that the form $F_2^3\omega^{\tens 2}$ has a zero
along~$E$, so that (over $U$, i.e., near~$P$) $F_2^3\omega^{\tens
2}\in H^0(2K-E)$. By assumption, $F_n\omega^{\tens n'}$ is a section
of $\psi^*(n'K_{\At})=n'(K-\alpha E)$. Therefore
$F_2^nF_n\omega^{\tens n}=(F_2^3\omega^{\tens 2})^{\tens
n'}\tens(F_n\omega^{\tens n'})$ gives a section of $n'(2K-E)+n'(K-\alpha
E) = nK-n'(1+\alpha)E$. But $E$ came from resolving a cyclic quotient
singularity (see [Br]) and cyclic quotient singularities are log
terminal (see, for instance, [CKM], Proposition~6.9): that is,
$\alpha>-1$. So we have a section of~$nK$, as required.~\qed

\noindent{\it Remark.\/} It is possible to show that $f(\bdw)$ has a
removable singularity along~$E$ by $L_2$ methods, following Freitag
([Fr]) and Sakai ([Sak]). $F_2$ is $L_3$-integrable and this implies
that $F_2\omega^{2/3}$ is $L_{2/3}$, that is, $\eta_2=F_2^3\omega^2$
satisfies $\left\{\int_U(\eta_2\wedge\bar\eta_2)\right\}<\infty$. The
argument of [Sak], Lemma~1.1, shows that $a_r=0$ if $r\le-k$ for an
$L_k$~form. In this case all we get is $a_r=0$ for $r<0$, which is not quite
good enough. I do not know whether a better bound can be obtained by a
refinement of the method of [Sak] in the case of non-integral~$k$.

\startsection 6 Conclusions.

We can now establish the main result by assembling the results we have
proved so far.

\thm 6.1 Any algebraic compactification of the moduli space $\Ap$ of
abelian surfaces with a polarisation of type $(1,p)$ for $p$ a prime
is of general type if~$p\ge 173$.

By Proposition~2.1 there is a weight~$2$ cusp form if $p>
71$. Suppose $n$ is sufficiently divisible. The obstructions to
extending $\eta=F_2^nF_n\omega^{\tens n}$ to the whole of a resolution
are the numerical obstruction given in Theorem~4.10 and the condition
in Theorem~5.3 that $F_n\omega^{\tens n'}$ (where $n=3n'$) should
extend to the boundary near an isolated non-canonical singular point
in~$E(h)$. This second condition can be met, according to
Proposition~3.1, by assuming that $F_n$ itself expressible as
$F_2^{n'}F_{n'}$, where $F_{n'}$ is a cusp form of weight~$n'$. So we
need a form $\eta=F_2^{4n'}F_{n'}\omega^{\tens 3n'}$. The dimension of
the space of such forms is ${{p^2+1}\over{8640}}n'^3+O(n'^2)$, by
Proposition~2.2. So, comparing
this with the obstruction in Theorem~4.10, we see that the moduli
space is of general type as long as $p>71$ and
$$
{{p^2+1}\over{8640}}>{{7}\over{2}}-{{9}\over{p}}
$$
which gives $p\ge 173$.~\qed

\noindent{\it Remarks.\/} a) This bound is not at all likely to be
sharp. There will in general be weight~$n$ forms $F_n$ not expressible
as $F_2^{n'}F_{n'}$ such that the corresponding differential form
$F_n\omega^{\tens n'}$ nevertheless extends to the general point of
the boundary. And for large~$p$ there will be many independent
weight~$2$ forms and therefore one expects there to be many more forms
expressible as $F_2^{n'}F_{n'}$ than we have actually written
down. However, it is not clear that such an expression should be
unique if it exists, so we cannot easily determine the dimension of
the relevant space.

\noindent{\phantom{\it Remarks.\/}} b) In [GH] and [GS] the use of the
weight~$2$ form was purely to get a better bound, but here we have
made essential use of the existence of such a form (in Theorem~5.3) in
order to get any bound at all.

It is probably true that the result of [Bor] still holds if
$\Sp(4,\ZZ)$ is replaced by $\Gp$. If this could be shown (it would be
enough to show it for the primes not covered by Theorem~6.1) then all
but finitely many moduli spaces of prime-polarised abelian surfaces
would be of general type. Practically, it would be best to improve
Theorem~6.1 first, so as to have only a few cases to deal with, if
possible. Ideally one would like to drop the restriction that $p$ be a
prime, but the singularities and the boundary become more complicated.

There ought to be a similar result even for $t$ composite, but the
singularities become complicated (there is a partial description in
[Br]): instead of four singular curves in $H_1$ we have to handle two
such curves for each ordered pair $(t',t'')\in\NN\times\NN$ such that
$t't''=t$. In general one expects ${\cal A}^*_t$ to to become closer
to general type as $t$ gets larger, but if $t$ has many prime factors
that tends to make ${\cal A}^*_t$ have lower Kodaira dimension. For
instance, it is shown in~[G] that $p_g({\cal A}^*_{13})\ge 1$ and
therefore $\kappa({\cal A}^*_{13})\ge 0$, but Gross and Popescu ([GP]) have
recently shown that ${\cal A}^*_{14}$ is unirational. In some special
cases, though, we can say a bit more.

\corollary 6.2 If $q\in\NN$ and $p\ge 173$ is a prime, then any
compactification of ${\cal A}^*_{pq^2}$ is of general type.

\pf If we put $Q=\diag (1,q^{-2},1,q^2)\in\Sp(4,\QQ)$ then
$Q\Gamma_{pq^2}Q^{-1}\contin\Gp$. This means that there is a
surjective morphism ${\cal A}_{pq^2}\to\Ap$ and this can be extended
to suitable smooth compactifications. Hence ${\cal A}^*_{pq^2}$ is of
general type as long as $\At$ is.~\qed

\startsection References\par
\exhyphenpenalty100
\frenchspacing
\item{[AS]}M.F.~Atiyah, I.M.~Singer, The index of elliptic
operators,~III, Ann.\ of\ Math.\ {\bf 87} (1968), 546--604.
\item{[Bor]}L.A.~Borisov, Finiteness theorem for $\Sp(4,\ZZ)$,
Preprint, Ann Arbor, Mich., 1994.
\item{[Br]}H.-J.~Brasch, Modulr\"aume abelscher Fl\"achen, Thesis,
Erlangen, 1994.
\item{[CKM]}H.~Clemens, J.~Koll\'ar, S.~Mori, {\sl Higher dimensional
complex geometry}, Ast\'erisque\ {\bf 166}, Soci\'et\'e Math\'ematique
de France, Paris, 1988.
\item{[D]}V.~I. Danilov, The Geometry of Toric Varieties, Russian\
Math.\ Surveys\ {\bf 33} (1978), 97--154.
\item{[EZ]}M.~Eichler, D.~Zagier,
{\sl The theory of Jacobi forms},
Progress in Math.\ {\bf 55}, Birkh\"auser, Boston, 1985.
\item{[Fr]}E.~Freitag,
{\sl Siegelsche Modulfunktionen}, Ergebnisse der mathematischen
Wissenschaften\ {\bf 254}, Springer, Berlin, 1983.
\item{[G]}V.~Gritsenko, Irrationality of the moduli spaces of polarized
abelian surfaces, Int.\ Math.\ Research\  Notices\ {\bf 6} (1994),
235--243.
\item{[GH]}V.~Gritsenko, K.~Hulek, Appendix to the paper
``Irrationality of the moduli space of polarized
abelian surfaces'', Preprint, Hannover 1994.
\item{[GP]}M.~Gross, S.~Popescu, in preparation.
\item{[Hir]}F.~Hirzebruch, Elliptische Differentialoperatoren auf
Mannigfaltigkeiten. Arbeitsgemeinschaft f\"ur Forschung des Landes
Nordrhein--Westfalen, Wissenschaftliche Abhandlungen Bd.\ {\bf 33},
583--608, Westdeutsche Verlag, K\"oln 1966. Also in: F.~Hirzebruch, Gesammelte
Abhandlungen Bd.~II, 24-49. Springer, Berlin 1987.
\item{[HKW1]}K.~Hulek, C.~Kahn, S.~Weintraub, Singularities of the
moduli spaces of certain abelian surfaces, Compositio\ Math.\ {\bf 79}
(1991), 231--253.
\item{[HKW2]}K.~Hulek, C.~Kahn, S.~Weintraub, {\sl Moduli spaces of
abelian surfaces: compactification, degeneration and theta functions.}
Expositions in Mathematics\ {\bf 12}, De Gruyter, Berlin, 1993.
\item{[HS]}K.~Hulek, G.~K. Sankaran, The Kodaira dimension of certain moduli
spaces of abelian surfaces,
Compositio\ Math.\ {\bf 90} (1994), 1--36.
\item{[HW]}K.~Hulek, S.~Weintraub, Bielliptic abelian surfaces,
Math.\ Ann.\ {\bf 283} (1989), 411--429.
\item{[LB]}H.~Lange, Ch.~Birkenhake,
{\sl Complex abelian varieties}, Grundlehren\ {\bf 302}, Springer,
Berlin, 1992.
\item{[Mor]}D.~Morrison, Canonical quotient singularities in dimension three,
Proc.\ A.M.S.\ {\bf 93} (1985), 393--396.
\item{[MorS]}D.~Morrison, G.~Stevens, Terminal quotient singularities
in dimensions three and four,
Proc.\ A.M.S.\ {\bf 90} (1984), 15--20.
\item{[Mum]}D.~Mumford, {\sl Abelian varieties}, Oxford University
Press, Bombay, 1970.
\item{[Od]}T.~Oda,
{\sl Convex bodies and algebraic geometry}, Ergebnisse der Mathematik
und ihrer Grenzgebiete, 3.~Folge, Bd\ {\bf 15}, Springer, Berlin, 1988.
\item{[O'G]}K.~O'Grady, On the Kodaira dimension of moduli spaces of
abelian surfaces, Compositio\ Math.\ {\bf 72} (1989), 121--163.
\item{[Sak]}F.~Sakai, Kodaira dimensions of complements of divisors.
In {\sl Complex analysis and algebraic geometry} (W.L.~Baily \&
T.~Shioda, eds.), 239--258, Cambridge University Press, Cambridge
1977.
\item{[SC]}A.~Ash, D.~Mumford, M.~Rapoport, Y.~Tai, {\sl Smooth
Compactification of Locally Symmetric Varieties}, Math.\ Sci.\ Press,
Brookline, Mass., 1975.
\item{[Se]}J.-P.~Serre, {\sl Repr\'esentations lin\'eaires des groupes
finis}, Hermann, Paris, 1967.
\item{[SZ]} N-P.~Skoruppa, D.~Zagier,
Jacobi forms and a certain space of modular forms,
Invent.\ Math.\
{\bf 94} (1988), 113--146.
\item{[T]}Y.~Tai, On the Kodaira dimension of the moduli spaces of
abelian varieties, Invent.\ Math.\ {\bf 68} (1982), 425--439.
\item{[YPG]}M.~Reid, Young person's guide to canonical singularities.
In {\sl Algebraic Geometry, Bowdoin, 1985} (D. Eisenbud, ed.), Proc.\
Symp.\ Pure\ Math.\ {\bf 46} (1987), 345--416
\item{[Z]}J.~Zintl, Thesis, Hannover, in preparation.
\raggedright
\medskip
\noindent G.~K. Sankaran, Department of Pure Mathematics and
Mathematical Statistics,\hfil\break
16,~Mill~Lane, Cambridge CB2~1SB, England.

\end